\documentclass[11pt]{article}
 
\usepackage[applemac]{inputenc} 

\usepackage{amsthm,amssymb,amsbsy,amsmath,amsfonts,amssymb,amscd,mathrsfs}
\usepackage[normalem]{ulem}
\usepackage{cite}
\usepackage{url}
\usepackage{tikz}

\long\def\symbolfootnote[#1]#2{\begingroup%
\def\thefootnote{\fnsymbol{footnote}}\footnote[#1]{#2}\endgroup} 

\usepackage{graphicx}
\DeclareGraphicsExtensions{.jpg}
\usepackage{color}
\usepackage{fullpage}
\newtheorem{theorem}{Theorem}
\newtheorem{corollary}[theorem]{Corollary}
\newtheorem{lemma}[theorem]{Lemma}
\newtheorem{proposition}[theorem]{Proposition}

\theoremstyle{definition} 
\newtheorem{definition}[theorem]{Definition}
\newtheorem{example}[theorem]{Example}
\theoremstyle{remark}
\newtheorem{remark}[theorem]{Remark}
\newcommand{\bt}{\begin{theorem}}
\newcommand{\et}{\end{theorem}}
\newcommand{\bl}{\begin{lemma}}
\newcommand{\el}{\end{lemma}}
\newcommand{\bp}{\begin{proposition}}
\newcommand{\ep}{\end{proposition}}
\newcommand{\bc}{\begin{corollary}}
\newcommand{\ec}{\end{corollary}}
\newcommand{\bdeff}{\begin{definition}}
\newcommand{\edeff}{\end{definition}}
\newcommand{\brem}{\begin{remark}}
\newcommand{\erem}{\end{remark}}
\newcommand{\bex}{\begin{example}}
\newcommand{\eex}{\end{example}}

\renewcommand{\phi}{\varphi}

\newcommand{\bi}{\begin{itemize}}
\newcommand{\iii}{\item}
\newcommand{\ei}{\end{itemize}}
\newcommand{\bd}{\begin{description}}
\newcommand{\ed}{\end{description}}

\newcommand{\bqn}{\begin{eqnarray}}
\newcommand{\eqn}{\end{eqnarray}}
\newcommand{\eqnn}{\nonumber\end{eqnarray}}

\newcommand{\nn}{\nonumber}
\newcommand{\ba}[1]{\begin{array}{#1}}
\newcommand{\ea}{\end{array}}
\newcommand{\wh}[1]{\widehat{#1}}

\newcommand{\al}{\alpha}
\newcommand{\eps}{\varepsilon}

\newcommand{\R}{\mathbb{R}}
\newcommand{\N}{\mathbb{N}}

\newcommand{\ex}{\exists\,}

\newcommand{\dive}{\text{div}}
\newcommand{\grad}{\nabla}

\def\dist{{d}}

\def\H{{\mathcal H}}
\def\1{{\mathbf 1}}
\def\ud{{\mathrm{\,d}}}
\def\vd{{\mathrm{d}}}

\def\cut{\mathrm{Cut}}

\def\hess{\mathrm{Hess}\,}
\def\vol{\,\mathrm{d vol}}
\def\thcut{\theta_{\mathrm{cut}}}
\def\tcut{t_{\mathrm{cut}}}

\usepackage{graphicx}
\DeclareGraphicsExtensions{.jpg}
\usepackage{color}
\usepackage{fullpage}
 \usepackage[usenames,dvipsnames]{pstricks}
 \usepackage{epsfig}
 \usepackage{pst-grad} 
 \usepackage{pst-plot} 

\begin{document}
\begin{center} \noindent
{\LARGE{\sl{\bf Small time heat kernel  asymptotics at the cut locus on surfaces of revolution}}}
\vskip 0.6 cm
Davide Barilari
\\ 
{\footnotesize CNRS, CMAP Ecole Polytechnique, Paris, 
 Equipe INRIA GECO Saclay-\^Ile-de-France, Paris, France  {\tt barilari@cmap.polytechnique.fr}}\\
\vskip 0.3cm
Jacek Jendrej
{\footnotesize \\CMAP Ecole Polytechnique, Paris,  France \\ {\tt jacek.jendrej@polytechnique.edu }}

\vskip 0.3cm
\end{center}

\vskip 0.1 cm
\begin{center}
\today
\end{center}
\vskip 0.2 cm
\begin{abstract} 
In this paper we investigate the small time heat kernel asymptotics on the cut locus on a class of surfaces of revolution, which are the simplest 2-dimensional Riemannian manifolds different from the sphere
with non trivial cut-conjugate locus.
We determine the degeneracy of the exponential map near a cut-conjugate point and present the consequences
of this result to the small time heat kernel asymptotics at this point.
These results give a first example where the minimal degeneration of the asymptotic expansion at the cut locus is attained.
\end{abstract}


\section{Introduction}

The study of the heat kernel and the geometric analysis of the Laplace-Beltrami operator on Riemannian manifolds is a very old problem
and it has been an object of an increasing attention during the last century.

In particular it has undergone a strong development in the last three decades during which many results have been obtained relating the properties of the heat kernel (such as Gaussian bounds, asymptotics, etc.) to the geometry of the manifold itself (optimality of geodesics, bounds on curvature, etc.). A comprehensive introduction to the subject can be found in \cite{rosenberg,grigoryan} and references therein.

In this spirit, the fundamental relation between the metric structure of the manifold and the properties of the heat kernel is provided by a result of Varhadan, stating that the main term in the small time asymptotic expansion of the heat kernel is given by the Riemannian distance. 
\bt[see \cite{behavior-of-heat-kernel}] Let $M$ be a complete Riemannian manifold,  $d$ be  the Riemannian distance and $p_{t}$ its heat kernel. For $x,y\in M$ we have 
 \bqn \nn
 \label{eq-varadh}
\lim_{t\to0}4 t \log p_t(x,y)=-d^2(x,y),
\eqn
uniformly for $(x,y)$ in every compact subset of $M\times M$.
\et
Besides this result, which gives the basic estimate for the heat kernel and holds for every pair of points on the manifold, one can study how the small time asymptotics of the heat kernel reflects the singularities of the Riemannian distance, to get refined results. By singularities of the distance here we mean the set of points where the Riemannian distance is not smooth. In particular we are interested in how the presence of the cut locus (the set of points where the geodesics lose global optimality) and the conjugate locus (the set of points where the differential of the exponential map is not invertible) can affect the asymptotic expansion of the heat kernel.

This relation has already been  established in the case of Riemannian geometry, as stated in \cite{molchanov,neel,davide-heat} (see also \cite{hsu} for a probabilistic approach to the Riemannian heat kernel). We recall here a corollary of these results that is contained in \cite{davide-heat} (where the results are stated in the more general framework of sub-Riemannian geometry \cite{montgomerybook,hypoelliptic}), saying that the leading term in the expansion of $p_{t}(x,y)$ for $t\to 0$ other than the exponential one depends on the structure of the minimal geodesics connecting $x$ and $y$.
\bt[see \cite{davide-heat}]
\label{corollarioMAIN0}
Let $M$ be a complete Riemannian manifold, let $d$ denote the Riemannian distance and $p_{t}$ the heat kernel. For $x,y\in M$ we have the bounds

\begin{equation}\nn
\frac{C_1}{t^{n/2}}  e^{-d^{2}(x,y)/4t} \leq p_t(x,y)
\leq\frac{C_2}{t^{n-(1/2)}}  e^{-\dist^{2}(x,y)/4t}, \quad \text{ for }
0<t<t_0.
\end{equation}
for some $t_{0}>0$, where $C_{1},C_{2}>0$ depend on
$M$, $x$, and $y$. Moreover\\[0.2cm]
(i) if $x$ and $y$ are not conjugate along any
minimal geodesic joining them, then
\[
p_t(x,y) = \frac{C+O(t)}{t^{n/2}}  e^{-\dist^{2}(x,y)/4t}, \qquad \text{
for } 0<t<t_0.
\]
(ii) if $x$ and $y$ are conjugate along at least one minimal geodesic
connecting them, then
\bqn \label{eq:bound}
p_t(x,y)\geq \frac{C}{t^{(n/2)+(1/4)}}  e^{-\dist^{2}(x,y)/4t}, \qquad
\text{ for } 0<t<t_0.
\eqn
\et

Note that (i) holds in particular for $y$ close enough to $x$ since the exponential map starting from a point is always a local diffeomorphism in a neighborhood of the point itself. 

The goal of this paper is to show some explicit examples of two dimensional manifolds (that include ellipsoids of revolution as a particular case) where the inequality (ii) in Corollary \ref{corollarioMAIN0} is sharp at some point $y$ that is both cut and conjugate with respect to $x$ (in the sense that the exponent $\alpha$ in the asymptotics $p_{t}(x,y)\sim t^{-\alpha}e^{-\dist^{2}(x,y)/4t}$ satisfies $\alpha=n/2+1/4$).

To this end, we analyze the small time asymptotics on some particular class of two dimensional surfaces
which are the simplest perturbations of the two-dimensional sphere.
Indeed the first non trivial example of a complete Riemannian manifold where an asymptotic expansion at the cut locus is computed is given by the standard two-dimensional sphere $S^{2}$. In this case the cut locus of a point $x$ coincides with the conjugate locus and collapses to the antipodal point $\wh{x}$ to $x$. The heat kernel on the two sphere was first computed in \cite{fischer2s} and an elementary computation  shows that the small time asymptotics at these points is given by
\bqn
\label{eq:s2}
p_t(x,\wh{x})\sim \frac{1}{t^{3/2}} e^{-\dist^{2}(x,\wh{x})/4t}, \qquad
\text{ for } t\to0.
\eqn
Since, to the authors' best knowledge, there is no explicit proof of this expansion (although it is stated in \cite{molchanov}) we wrote it in  Appendix \ref{app:s2}.
Notice that in this case $n=2$, thus the exponent in the expansion $p_{t}(x,y)\sim t^{-\alpha}e^{-\dist^{2}(x,y)/4t}$ is not the minimal one (cf. \eqref{eq:bound}), as a consequence of the fact that the cut-conjugate point is reached by a one parametric family of optimal geodesics. (See also Remark \ref{rem:sphere}.)

The next class of examples that is natural to consider is the one of two-dimensional surfaces of revolution, in which of course  ellipsoids of revolution are the main model. 

The determination of the cut and conjugate loci on a complete surface, even a two-dimensional surface of revolution, is a classical but difficult problem in Riemannian geometry. Even on an ellipsoid of revolution, the computation is not a standard exercise. In \cite{bergerbook}, the foreseen conjugate and cut loci were given as a conjecture and the first proof appeared only recently in \cite{itoh-kiyohara}. For other results about the conjugate and cut locus on surfaces of revolution (and some generalizations) one can see also \cite{sinclair-tanaka,bonnardcaillausinclair,itohcut2}.

On an oblate ellipsoid of revolution the cut locus of a point different from the pole is a subarc of the antipodal parallel.
On a prolate ellipsoid it is a subarc of the opposite meridian. In the first case the Gaussian curvature is monotone increasing from the north pole to the equator and decreasing in the second case.

This result is also a particular case of a more general one contained in \cite{sinclair-tanaka} about the cut locus of a so-called \emph{two-sphere of revolution}.

\bt[see  \cite{sinclair-tanaka}] Given a smooth metric on $S^{2}$ of the form $dr^{2} +m^{2}(r)d\theta^{2}$, where $r$ is the distance along the meridian and $\theta$ the angle of revolution and $m$ is smooth, assume the following:
\bi
\iii[(i)] $m(2a-r)=m(r)$ (i.e. reflective symmetry with respect to the equator, where $2a$ is the distance between poles).
\iii[(ii)] the Gaussian curvature is monotone non-decreasing (resp. non-increasing) along a meridian from the north pole
to the equator.
\ei
Then the cut locus of a point different from the pole is a simple branch located on the antipodal parallel (resp. opposite meridian).
\et


In this paper we will deal only with the cut locus of a point that belongs to the equator of the two sphere of revolution. Due to this restriction we can require the reflective symmetry with respect to the equator only locally, that means that our result is stated under the following assumptions:

\begin{enumerate}
\item[(A1)] $m(r)=\psi((a-r)^2)$  in some neighborhood of $r=a$, where $\psi$ is a smooth function. \label{ass:reflection}
\item[(A2)] the Gaussian curvature is monotone non-decreasing along a meridian from the pole
to the equator.\label{ass:gaussian}
\end{enumerate}
Notice that in the condition (A2) we specify the Gaussian curvature to be monotone increasing. For instance, in the class of ellipsoid of revolution, this means that we are considering only the oblate ones.

To prove our result about the asymptotics of the heat kernel, we also need a \emph{nonsingularity} assumption of the  two-sphere, that is the following:
\bi
\iii[(A3)] $K''(a)\neq 0$ where $K(r_0)$ is the Gaussian curvature (that is constant) along the parallel $r = r_0$.
\ei
Notice that an oblate ellipsoid of revolution is a typical example of a two-sphere of revolution satisfying all the three assumptions (A1), (A2) and (A3). (See also Section \ref{sec:ellipsoid}.)

Using these results, we analyzed the structure of the degeneracy of the exponential map on such a  two-sphere of revolution around a cut-conjugate point on the equator proving that, under these assumptions, the following asymptotic expansion holds.

\bt\label{thm:main}\label{thm:main2}
Let $M$ be a complete two-sphere of revolution satisfying assumptions (A1)-(A3) above, and let $d$ be its Riemannian distance and $p_{t}$ its heat kernel. Fix $x\in M$ along the equator and let $y$ be a cut-conjugate point with respect to $x$. Then we have the asymptotic expansion
$$
p_t(x,y)\sim \frac{1}{t^{5/4}} e^{-\dist^{2}(x,y)/4t}, \qquad
\text{ for } t\to0.
$$
\et

This result is a consequence of the fact that the degeneracy of the exponential map is in relation with the asymptotics of the heat kernel via the asymptotic expansion of the so-called \emph{hinged energy function} (see Section \ref{s:heat} for the precise definition and \cite[Theorem 24]{davide-heat} for the aforementioned relation) which in the two dimensional case can be explicitly analyzed.

\subsection{Structure of the paper}
In Section \ref{s:prel} we introduce some preliminary material. In particular Section \ref{s:heat} contains the results about the relation between the heat kernel asymptotics and the degeneracy of the exponential map, while Section \ref{s:2s} recalls the formal definition and the basic properties of the two-spheres of revolution. Sections \ref{s:main1} and \ref{s:main2} contain the proof of the main result of the paper. In Section \ref{sec:ellipsoid} these results are expressed in terms of the geometry of ellipsoids.
 Finally in Appendix \ref{app:s2}, given the formula for the heat kernel on $S^{2}$, we compute its expansion at the cut locus, while in Appendix \ref{a:lemmas} we prove some technical lemmas that are needed in the proof of the main result.

\section{Preliminary material}\label{s:prel}
In this section we introduce some preliminary material that is needed for the proof of our main theorem. In particular we  recall the main results about the relation between the heat kernel asymptotics and the degeneracy of the exponential map. For more details one can see \cite{davide-heat}. A description of these results in the Riemannian case can be also found in\cite{molchanov, neel}.

\subsection{Hinged energy function and cut-conjugate points}\label{s:heat}

In what follows, wherever not specified, $M$ is an $n$-dimensional complete Riemannian manifold and $d$ denotes the Riemannian distance.\\

{\bf Notation.} If $x\in M$ and $v\in T_{x}M$ we denote by $\gamma_{v}(t)$ the geodesic starting from $x$ with initial velocity $v$. 
Moreover we denote by $\exp:TM\to M$ the exponential map defined by $\exp(v) = \gamma_ v(1)$.
If the initial point is fixed, we use the notation $\exp_x = \exp\vert_{T_xM}$.
 \bdeff
Let $x\in M$. The \emph{energy function} from $x$ is the function $E_x: M\to\R$ defined by 
\begin{equation*}\label{eq:en-function}
E_x(y) := \frac{1}{2}\dist(x, y)^2.
\end{equation*}
\edeff
The gradient of the energy function at a point $y\in M$, at a differentiability point, is easily computed by means of the unique geodesic joining $x$ and $y$. We have the following proposition (see  \cite{gallot-geometry}):
\bp\label{prop:en-deriv}
The function $E_{x}$ is smooth in the set $V_{x}\subset M$ of points that are reached from $x$ by a unique minimizing and non conjugate geodesic. Moreover, if $y=\exp_{x} (v)\in V_x$ then $\grad E_x(y) = \dot\gamma_v(1) \in T_y M$.
\ep
%

Now we define the main object that is involved in the asymptotics of the heat kernel.
\bdeff
Let $x, y\in M$. The \emph{hinged energy function} $h_{x, y}:M\to\R$ is defined by
\begin{equation*}
h_{x, y}(z) := E_x(z) + E_y(z), \qquad z\in M.
\end{equation*}
\edeff
Observe that the function $h_{x, y}$ is smooth on $V_x\cap V_y$. In the next proposition we characterize its minima.
\bp\label{prop:hinged-minimum}
Let $x, y\in M$ and let $\gamma:[0, 2]\to M$ be a minimal geodesic from $\gamma(0) = x$ to $\gamma(2) = y$.
Then $\min h_{x, y} = \frac{1}{4}\dist(x, y)^2$ and this minimum is attained at the point $z_0 := \gamma(1)$.
Moreover every global minimum of $h_{x, y}$ is a midpoint of some minimal geodesic joining $x$ and $y$.
\ep
\begin{proof}
Let $z\in M$, and set $a := \dist(x, z_0)$, $b := \dist(z_0, y)$, $\alpha := \dist(x, z)$, $\beta := \dist(z, y)$.
Then $a = b = \frac{1}{2}\dist(x, y)$, so $h_{x, y}(z_0) = \frac{1}{2}a^2+\frac{1}{2}b^2 = \frac{1}{4}\dist(x, y)^2$.
From the triangle inequality we have $\alpha + \beta \geq \dist(x, y)$, thus
\begin{equation}\label{eq:111}
h_{x, y}(z) = \frac{\alpha^2+\beta^2}{2} \geq \left(\frac{\alpha+\beta}{2}\right)^2 \geq \frac{1}{4}\dist(x, y)^2.
\end{equation}
Suppose now that the equality in \eqref{eq:111} holds. Then $\frac{\alpha^2+\beta^2}{2} = \left(\frac{\alpha+\beta}{2}\right)^2$ and $\alpha + \beta = \dist(x, y)$. From this we deduce that $z$ lies on a minimal curve from $x$ to $y$ and
$\alpha = \beta$.
Then this curve is a geodesic and $z$ is its midpoint.
\end{proof}
\begin{remark}
Notice that the midpoint of a minimal geodesic lies in the ``good'' set $V_x\cap V_y$ even if $y\in \cut(x)$.

\end{remark}

The next proposition shows how the degeneracy of $\exp_x$ near $y$ is reflected by the behavior of $h_{x,y}$ near the midpoint
of a geodesic joining them, which is assumed to lie in the ``good'' region.
\bp\label{prop:transport}
Fix $x, y\in M$ and let $\alpha:(-\epsilon, \epsilon)\to V_x\cap V_y\subset M$ be a smooth curve such that
$\alpha(0)$ is a critical point of $h_{x,y}$. Assume moreover that $\ex k\in\N$ such that
$\grad h_{x, y}(\alpha(s)) = s^{k}u(s)$,
where $u$ is a smooth vector field along the curve $\alpha$ with $u(0) \neq 0$.

Let $\gamma_{v(s)}$ be the minimal geodesic joining $x$ to  $\alpha(s)$ in time 1 and set $\beta(s) := \gamma_{v(s)}(2)$.\\
Then $\beta:(-\epsilon, \epsilon)\to M$ is a smooth function, $\beta(0) = y$ and
\begin{equation}\label{eq:xxx}
\lim_{s\to 0}\frac{\dot\beta(s)}{s^{k-1}} = k\ud\exp \left(\dot{\gamma}_{v(0)}(1)\right)u(0),
\end{equation}
where in the last formula we identify $T_{\alpha(0)}M$ with $T_{\dot{\gamma}_{v(0)}(1)}(T_{\alpha(0)}M)$.
\ep
\begin{remark}\label{rem:v-smooth}
Observe that under our assumptions $v(s)$ is a smooth curve $v:(-\epsilon, \epsilon)\to T_x M$
and every such curve $v(s)$ gives rise to the corresponding curve $\alpha(s) := \gamma_{v(s)}(1)$.
\end{remark}
\brem If $\alpha(s)$ is the curve of midpoints of a one parameter family of optimal geodesics joining $x$ to $y$ then $\dot \beta(s)\equiv 0_y$ since the curve $\beta(s)$ is constant.
\erem
\begin{proof}
The first part is obvious, since $\beta$ is a composition of smooth functions. Let $\gamma_{w(s)}:\R\to M$ be the optimal geodesic from $\gamma_{w(s)}(0) = y$ to $\gamma_{w(s)}(1) = \alpha(s)$.
It follows from Proposition \ref{prop:en-deriv} that
\begin{equation}\label{eq:lambda}
\dot{\gamma}_{v(s)}(1) = -\dot{\gamma}_{w(s)}(1) + \grad h_{x, y}(\alpha(s))= -\dot{\gamma}_{w(s)}(1) + s^k u(s).
\end{equation}
We want to use the Taylor expansion of $\exp$ near $\dot{\gamma}_{v(0)}(1) = -\dot{\gamma}_{w(0)}(1)\in TM$.

To this end we introduce on $T M$ local canonical coordinates $(z, v) \in \R^{2n}$
in a neighborhood of $0_{\alpha(0)} \in TM$.
In these coordinates we write
\begin{align*}
\dot{\gamma}_{v(0)}(1) = -\dot{\gamma}_{w(0)}(1) &= (0, c), \\
\dot{\gamma}_{v(s)}(1) &=(z(s), c+ a(s)), \\
-\dot{\gamma}_{w(s)}(1) &= (z(s), c+ b(s)), \\
u(s) &= (z(s), u(s)),
\end{align*}
so that $a(s) - b(s) = s^k u(s)$.
Let
\begin{align*}
\exp_i(z, c+h) &= \left(\vd\exp_{(0, c)}(z, h)\right)_i + \phi_2\left(z, h\right)^2 + \dots + \phi_{k}\left(z, h\right)^k+O(|(z, h)|^{k+1})
\end{align*}
be the Taylor expansion of the $i$-th coordinate of $\exp$ near $\dot{\gamma}_{v(0)}(1) = (0, c)$ (here $\phi_j$ is some $j$-linear map $(\R^{2n})^j \to \R$ for $j = 2, \ldots, k$).
%

Observe that $\exp_i\left(z(s), c+b(s)\right) = y_i$ and
$\exp_i\left(z(s), c+a(s)\right) = \beta_i(s)$.
Also, for $j = 2, \ldots, k$ we have
\begin{equation*}
\left|\phi_{j}\left(z(s), a(s)\right)^j - \phi_{j}\left(z(s), b(s)\right)^j\right| = O(s^{k+1}),
\end{equation*}
because $\left|\left(z(s), a(s)\right) - \left(z(s), b(s)\right)\right| = \left|s^k\left(0, u(s)\right)\right|= O(s^k)$ and
$\left|\left(z(s), a(s)\right)\right| + \left|\left(z(s), b(s)\right)\right| = O(s)$.
Hence
\begin{equation*}
\beta_i(s) = y_i+\left(\vd\exp_{(0, c)}(0, a(s)-b(s))\right)_i + O(s^{k+1})=y_i+s^k\left(\vd\exp_{(0, c)}(0, u(s))\right)_i + O(s^{k+1}),
\end{equation*}
so that
\begin{align*}
\dot \beta(s) &= k s^{k-1}\vd\exp_{(0, c)}(0, u(s)) + O(s^k) \\
&= k s^{k-1}\vd\exp_{(0, c)}(0, u(0)) + O(s^k) \\
&= k s^{k-1}\vd\exp_{\dot{\gamma}_{v(0)}(1)}(u(0)) + O(s^k),
\end{align*}
and the conclusion follows.
\end{proof}

\bc\label{cor:conjugate}
Let $x, y\in M$ and let $z_0\in M$ be a critical point of $h_{x, y}$. Assume additionally that $z_0\in V_x\cap V_y$.
Then $y$ is conjugate to $x$ along some geodesic if and only if $z_0$ is a degenerate critical point of $h_{x, y}$.
\ec
\begin{proof}
Suppose first that $y$ is conjugate to $x$ along $\gamma_{v(0)}$. Then there exists a smooth curve
$v:(-\epsilon, \epsilon)\to T_x M$ such that $\dot\beta(0) = 0$. Let $\alpha:(-\epsilon, \epsilon)\to M$
be the corresponding curve of midpoints (see Remark \ref{rem:v-smooth}). Then in \eqref{eq:xxx}
we must have $k \geq 2$.
Thus $\hess h_{x, y}(\alpha(0))(v, \dot\alpha(0)) = 0$ for any $v$.

Conversely, if $\hess h_{x, y}(z_0)(v, w) = 0$ for all $v$, it suffices to take an arbitrary regular curve
$\alpha:(-\epsilon, \epsilon)\to V_x\cap V_y$ such that $\alpha(0)=z_0$ and $\dot\alpha(0)=w$
to obtain $\dot\beta(0) = 0$.
\end{proof}
\begin{remark}\label{rem:non-degeneracy}  $\hess h_{x, y}(z_{0})$ is never degenerate along the direction of a minimal geodesic connecting $x$ and $y$ (where $z_{0}$ is its midpoint). Indeed if $\gamma:[0,2]\to M$ is such a minimal geodesic, $z_{0}$ is its midpoint and $\alpha(s) := \gamma(1+s)$, then $h_{x, y}(\al(s)) = \frac{1}{4}\dist(x, y)^2(1+s^2)$. 
\end{remark}


Let $x\in M$ and let $\gamma:\R\to M$ be a geodesic with $\gamma(0) = x$. Suppose that $\gamma(2) = y$
is a conjugate point of $x$ along $\gamma$ and that $\gamma|_{[0, 2]}$ is minimal (such a point
is called a \emph{cut-conjugate point}).\\[-0.4cm]

It follows from Proposition \ref{prop:hinged-minimum} and Corollary \ref{cor:conjugate} that $z_0 := \gamma(1)$ is
a global minimum of $h_{x, y}$ and a degenerate critical point. Assume that $M$ has dimension 2.
Then, according to Remark \ref{rem:non-degeneracy}, the degeneracy is only in one direction,
which allows us to use the following result, called  Splitting Lemma or  Refined Morse Lemma
(it is a special case of \cite[Lemma 1]{gromoll-meyer}).
\bl\label{l:l}
Let $f: U\subset\R^n\to \R$ be defined in an open neighborhood of $0$ and assume that 0 is an isolated minimum of $f$ such that $f(0)=df(0)=0$ and $\dim \ker d^{2}f(0)=1$.
Then there exist coordinates $(x_1, \ldots,x_{n-1}, x_{n})$ around 0 such that
$f(x_1, \ldots,x_{n-1}, x_n) = x_1^2+\dots+x_{n-1}^2+g(x_n)$,
where $g:\R\to\R_+$ is smooth and $g(z) = O(z^4)$.
\el
\begin{remark}
The function $g$ is not unique, but its order of vanishing at $z=0$ is -- this is the maximal order of vanishing at $z=0$
of functions $s\mapsto f(\alpha(s)) - f(z_0)$ for a smooth curve $s\mapsto \alpha(s)$ such that $\alpha(0) = z_0$. 
\end{remark}
Applying Lemma \ref{l:l} to $h_{x, y}$, we obtain that there exists a smooth local coordinate system $(z_1, z_2)$ near $z_0$
such that $h_{x, y}(z_1, z_2) = h_{x, y}(z_0) + z_1^2 + g(z_2)$,
where $g:\R\to\R_+$ is a smooth function and $g(z) = O(z^4)$. Assume that $g$ vanishes at $z = 0$ with finite order $k+1\in\N$ (thus $k\geq 3$).

We define in these coordinates the smooth curve $\alpha(s) := (0, s)$. It is immediate that
$\grad h_{x, y}(\alpha(s))$\  $= s^k u(s)$ for some smooth vector field $u$ with $u(0)\neq 0$.
Conversely, there are no smooth curves $\alpha(s)$ such that $\grad h_{x, y}(\alpha(s)) = s^{k+1} u(s)$
with a smooth vector field $u$.
Combining this with Proposition \ref{prop:transport} and Remark \ref{rem:v-smooth} we obtain the following result.
\bc\label{cor:forcing-degeneracy}
Assume that $M$ has dimension 2. Let $x, y\in M$ and assume that $y$ is a cut-conjugate point of $x$ along $\gamma:[0,2]\to M$.
Then there exists a smooth curve $v: (-\epsilon, \epsilon)\to T_x M$ such that $\gamma_{v(0)} = \gamma$ and
\begin{equation*}
\dist(\gamma_{v(s)}(2), y) = O(s^k),
\end{equation*}
where $k+1$ is the order of vanishing at $0$ of the function $g$ described above.
It is impossible to obtain a higher order of vanishing.
\ec
This is the aforementioned relationship between the hinged energy function and the order of degeneracy of the exponential map.
We will show that for two-spheres of revolution satisfying assumptions (A1)-(A3), the order of degeneracy is lowest possible, that is $k+1 = 4$.\\

We  now briefly describe how the order of vanishing of the function $g$ is related to the small-time
asymptotics of the heat equation on $M$.

Let $M$ be a complete orientable $n$-dimensional Riemannian manifold.
We denote $\vol$ the volume form associated with the Riemannian metric on $M$ and compatible with the orientation
(which means it equals $1$ on a positively oriented orthonormal frame).

The \emph{Laplace-Beltrami operator} on $M$ is defined\footnote{recall that the \emph{divergence} $\dive X$ of a vector field $X$ is the unique function satisfying
$
L_X\vol = (\dive X)\vol
$, $L_{X}$ being the Lie derivative.} as $\Delta f := \dive\left(\grad f\right)$, $f\in C^{\infty}(M)$.
The \emph{heat operator} $e^{t\Delta}: L^2(M)\to L^2(M)$ is defined by the formula
\begin{equation*}
e^{t\Delta}f(p) = \int_M p_t(x, y)f(y)\vol(y),\qquad f\in L^{2}(M),
\end{equation*}
where $p_t(x, y)\in C^{\infty}(\R^{+}\times M\times M)$ is the \emph{heat kernel}, i.e. the unique smooth function  satisfying the following conditions.
\begin{align}
(\partial_t-\Delta_x)p_t(x, y) &= 0, \label{eq:heat-equation}\\
\lim_{t\to 0}\int_M p_t(x, y)f(y)\vol(y) &= f(x),\qquad \forall {f\in C^\infty(M)}\label{eq:heat-limit}.
\end{align}
Recall that the heat kernel exists and is unique. Moreover it satisfies $p_t(x, y) = p_t(y, x)>0$.

%
%
%

It is well known that if $M$ is compact then the operators $\Delta$ and $e^{t\Delta}$
are simultaneously diagonalizable and the proper vectors are smooth functions.
Their spectra determine the long-time behavior of the heat flow
and contain topological information (see for instance \cite{grigoryan,rosenberg,spectrebook}).

Here we investigate the short-time behavior, which turns out to be connected to the geometry of the manifold.
\bt\label{thm:davide}
Let $x, y\in M$ and assume that there exists only one minimal geodesic joining $x$ and $y$.
Let $z_0$ be the midpoint of this geodesic and assume that there exists a coordinate system $(z_1, \ldots, z_n)$ near $x$ such that
\begin{equation*}
h_{x, y}(z) = h_{x, y}(z_0) + z_1^{2m_1}+\dots+z_n^{2m_n}+o(|z_1|^{2m_1}+\dots+|z_n|^{2m_n}),
\end{equation*}
for some integers $1\leq m_1\leq\dots\leq m_n$.
Then there exists a constant $C>0$ (depending on $M$, $x$ and $y$) such that
\begin{equation*}
p_t(x, y) = \frac{C+o(1)}{t^{n-\sum_i\frac{1}{2m_i}}}\exp\left(-\frac{\dist(x, y)^2}{4t}\right).
\end{equation*}
\et
A full proof can be found in \cite{neel}, whereas the main ideas are already in \cite{molchanov}.
The result was extended to sub-Riemannian manifolds in \cite{davide-heat}. For other results about the small time asymptotics of the elliptic (and hypoelliptic) heat kernel one can see also \cite{benarous,benarousdiag,behavior-of-heat-kernel,leandre}.
\begin{remark}\label{rem:sphere}
As explained in \cite[Remark 2]{davide-heat}, if there exists a one-parameter family of minimal geodesics joining $x$ and $y$,
the theorem is still valid, but it should be understood that some $m_i$ is infinite.
For example let $M = S^2$ be a sphere of radius 1 and let $x, \wh{x}$ be two opposite poles. Then we have $m_1=1, m_2=\infty$,
so the asymptotics is (see Appendix \ref{app:s2})
\begin{equation*}
p_t(x, \wh{x}) = \frac{C+o(1)}{t^{3/2}}e^{-{\pi^2}/{4t}}.
\end{equation*}
\end{remark}
In particular, let us return to the situation of Corollary \ref{cor:forcing-degeneracy}. We obtain the following result.
\bc\label{cor:the-exponent}
Let $M$ be a $2$-dimensional orientable Riemannian manifold.
Let $x, y\in M$ and assume that there exists only one minimal geodesic joining $x$ and $y$.
Assume further that there exists no smooth curve $v:(-\epsilon, \epsilon)\to T_p M$ such that
$\gamma_{v(0)}:[0,2]\to M$ is the minimal geodesic from $x$ to $y$ and
\begin{equation*}
\dist(\gamma_{v(s)}(2), y) = o(s^3).
\end{equation*}
Then there exists a constant $C>0$ (depending on $M$, $x$ and $y$) such that
\begin{equation*}
p_t(x, y) = \frac{C+o(1)}{t^{5/4}}\exp\left(-\frac{\dist(x, y)^2}{4t}\right).
\end{equation*}
\ec

\subsection{Two-spheres of revolution} \label{s:2s}
In this section we introduce the class of two dimensional manifolds that we are going to study
and we recall their properties which are used in the sequel.
%
%
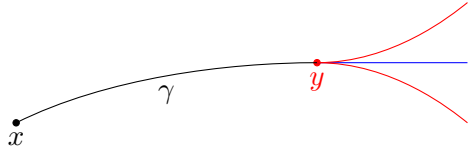
\begin{figure}
\begin{center}
\begin{tikzpicture}[scale=4]
\filldraw [black] (0,0) node[below] {$x$} circle (0.3pt);
\filldraw [red] (1,0.2) node[below] {$y$} circle (0.3pt);
\draw (0,0) .. controls (0.3,0.15) and (0.7,0.20) ..
node[below] {$\gamma$} (1,0.2);
\draw [blue] (1,0.2) -- +(0.5,0);
\draw [red] (1, 0.2) parabola +(0.5,0.2);
\draw [red] (1, 0.2) parabola +(0.5,-0.2);
\end{tikzpicture}
\caption{\label{fig:cut-conj}Picture of the cut locus (blue) and conjugate locus (red) from $x$ on a surface.}
\end{center}
\end{figure}
%
%

\begin{figure}
\begin{center}
\begin{tikzpicture}[scale=4]
\draw (0,0) -- node[below] {$b$} (1,0);
\draw (0,0) arc (180:90:1 and 0.5) -- (1,0);
\draw (1,0.5) node[above] {$p$}
arc (90:60:1 and 0.5) node[above] {$r$}
arc (60:30:1 and 0.5) node[right] {$(r, \theta)$}
-- node[below] {$m(r)$} (1,0.25);
\draw (2,0) arc (0:30:1 and 0.5);
\end{tikzpicture}
\caption{\label{fig:ellipse}Notational convention for a $2$-sphere of revolution.}
\end{center}
\end{figure}
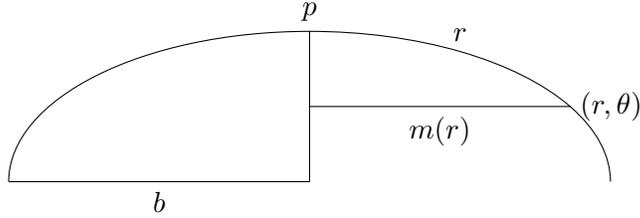

\bdeff
Let $M$ be a compact Riemannian manifold homeomorphic to $S^2$.
$M$ is called a \emph{$2$-sphere of revolution} if there
exists a point $p\in M$ called a \emph{pole} such that for any $q_1, q_2\in M$ satisfying $\dist(p, q_1)=\dist(p, q_2)$
there exists an isometry $f: M\to M$ satisfying $f(q_1)=q_2, f(p)=p$.
\edeff
Let $M$ be a $2$-sphere of revolution and let $p$ be a pole.
It can be proved that $p$ has a unique cut point $q$ which is also a pole \cite[Lemma 2.1]{sinclair-tanaka}. 
Therefore $M\setminus\{p, q\}$ can be parametrized by geodesic polar coordinates around $p$, which we denote $(r, \theta)$.
Let $M':=M\setminus\{p, q\}$.
This allows to express the Riemannian metric on $M'$ as $\ud s^2=\ud r^2+m(r)^2\ud \theta^2$ where $m$ is a positive function satisfying $\lim_{r\to 0}m(r) = 0$
\cite[p. 287, Proposition 3]{do-carmo-curves}.
 Geometrically, $m$  is the distance from the rotational axis (see also Figure \ref{fig:ellipse},
which shows a section of $M$ by a plane containing the rotational axis).

We denote $a := \frac{1}{2}\dist(p, q)$ and $b:=m(a)$. The set $\{(r, \theta)\in M: r=a\}$ is called the \emph{equator}.
Each set $\{(r, \theta)\in M: r = \mathrm{const}\}$ is called a \emph{parallel}, while
each set $\{(r, \theta)\in M: \theta=\mathrm{const}\}$ is called a \emph{meridian}.

Next, we recall a standard result on geodesics on surfaces of revolution.
\bp[see \cite{do-carmo-curves}]
Let $\gamma(t) = (r(t), \theta(t))$ be a unit speed geodesic on $M'$.
There exists a constant $\nu$, called the \emph{Clairaut constant} of $\gamma$, such that
\begin{equation}\label{eq:clairaut}
m(r(t))^2\dot\theta(t) = \nu.
\end{equation}
\ep
\begin{remark}
Observe that $m(r(t))\dot\theta(t) = \cos \eta(t)$, where $\eta(t)$ is the angle
between $\dot{\gamma}(t)$ and $\frac{\partial}{\partial\theta}\big|_{\gamma(t)}$. The constant $\nu$ can be interpreted as the angular momentum of the geodesic.
\end{remark}
We will study geodesics emanating from a point $x$ on the equator. Without loss of generality we assume that $\theta(x) = 0$.
\bp
Let $\nu\in\R$. Then we have the following properties:
\begin{enumerate}
\item{If $|\nu| > b$, then no geodesic on $M'$ emanating from $x$ satisfies \eqref{eq:clairaut}.}
\item{If $|\nu| = b$, then there is exactly one geodesic $\gamma: \R\to M'$
emanating from $x$ satisfying \eqref{eq:clairaut}. Its image is the equator.}
\item{If $0<|\nu| < b$, then there is exactly one geodesic $\gamma: \R\to M'$
emanating from $x$ satisfying \eqref{eq:clairaut} and $r'(0) < 0$.}
\end{enumerate}
\ep
\noindent For $\nu\in\R$ this unique geodesic will be denoted by $c_\nu$.
\begin{proof}
We assumed $\gamma(t)$ to be unit speed, that is $\dot r(t)^2+m(r(t))^2\dot\theta(t)^2=1$ for all $t$.
Thus \eqref{eq:clairaut} is equivalent to
\begin{equation}\label{eq:deriv-r}
\dot r(t)^2 = \frac{m(r(t))^2-\nu^2}{m(r(t))^2},
\end{equation}
which for $t = 0$ gives
\begin{equation}\label{eq:forr}
\dot r(0)^2 = \frac{b^2-\nu^2}{b^2}.
\end{equation}
\begin{enumerate}
\item{If $|\nu|>b$, we clearly get a contradiction.}
\item{If $|\nu| = b$, we get $\dot r(0) = 0$ and $\dot\theta(0) = \frac{\nu}{b^2}$.
There is a unique geodesic with initial tangent vector $(\dot r(0), \dot\theta(0)) = \left(0, \frac{\nu}{b^2}\right)$.
By symmetry this geodesic will not quit the equator. By \eqref{eq:clairaut} $\dot\theta(t)$ is constant,
so the geodesic covers the whole equator.}
\item{If $|\nu| < b$, we obtain $\dot r(0) = -\sqrt{\frac{b^2-\nu^2}{b^2}}$ (because we chose $\dot r(0)$ to be negative)
and $\dot\theta(0) = \frac{\nu}{b^2}$.
This determines the unique geodesic with Clairaut constant $\nu$. We have $|m(r(t))\dot\theta(t)| = |\cos\eta(t)|\leq 1$,
so $m(r(t)) \geq |\nu| > 0$. This ensures that the geodesic does not meet the poles.}
\end{enumerate}
\end{proof}

\noindent
{\bf Assumptions.} In what follows we make the following assumptions:
\begin{enumerate}
\item[(A1)] $m(r)=\psi((a-r)^2)$  in some neighborhood of $r=a$, where $\psi$ is a smooth function. \label{ass:reflection}
\item[(A2)] the Gaussian curvature is monotone non-decreasing along a meridian from the pole
to the equator.\label{ass:gaussian}
\iii[(A3)] $K''(a)\neq 0$ where $K(r_0)$ is the Gaussian curvature (that is constant) along the parallel $r = r_0$.
\end{enumerate}

The assumption (A1) means that the metric is invariant by reflection with respect to the equator, in a neighborhood of the equator itself. 

From  (A2) and the Gauss-Bonnet Theorem
it follows that the Gaussian curvature on the equator is strictly positive.
Thus from Lemmas 2.2 and 2.3 in \cite{sinclair-tanaka} one can obtain that $m$ is strictly increasing on $(0,a]$.
For every $\nu$ such that $0<|\nu|\leq b$, we denote by $R=R(\nu)$ the unique $R\in(0,a]$ such that $m(R)=|\nu|$.
It is the minimal geodesic distance of $c_\nu(t)$ from the pole $p$.

Theorem 4.1 in \cite{sinclair-tanaka} states that the cut locus of $x$ is a subset of the equator.
For $0<|\nu|<a$ the cut point along $c_\nu$ is the first point of intersection of $c_\nu$ with the equator.
This point is given by the formula \cite[p. 385]{sinclair-tanaka} $(r, \theta) = (a, \varphi(\nu))$, where
\begin{equation}\label{eq:function-varphi}
\varphi(\nu) := 2\int_R^a\frac{\nu\ud r}{m(r)\sqrt{m(r)^2-\nu^2}}.
\end{equation}

From now on, we will assume that $\nu \geq 0$ (the case $\nu \leq 0$ is symmetric).
It can be shown that $\varphi(\nu)$ is non-increasing for $\nu\in(0, b)$ \cite[Lemma 4.2]{sinclair-tanaka}.
Thus the cut-conjugate point of $x$ has coordinates $(r, \theta) = (a, \lim_{\nu\to b^-}\varphi(\nu))$.
We note $\thcut := \lim_{\nu\to b^-}\varphi(\nu)$ and $\tcut := b\thcut$, so that $c_b(\tcut)$ is the cut-conjugate point.

\begin{remark}
The geodesic $c_\nu$ starting from $x$ is determined by each of the parameters $\nu$, $R$ or $\eta$.
We recall the relationships between these parameters:  $\nu = m(R)$ and
$\cos\eta = \frac{\nu}{b}$. In particular $b-\nu$ is of order $\eta^2$ as $\eta\to 0$.
\end{remark}

\section{Proof of Theorem \ref{thm:main}}\label{s:main}
In what follows, we prove Theorem \ref{thm:main2} in two steps.  First we compute the asymptotic expansion of the function $\phi$ and then we apply this result to investigate the geodesic variations of the equator.

\subsection{Expansion of $\phi(\nu)$}\label{s:main1}

We will now derive the asymptotic expansion of $\varphi(\nu)$ for $\nu$ close to $b$
from the asymptotic expansion of $m(r)$ near $r=a$.
Let \begin{equation*}
\psi(s) = b-\alpha s+\beta s^2+O(s^3),
\end{equation*}
be the Taylor expansion of $\psi$, so that
\begin{equation*}
m(r) = \psi((a-r)^2)= b-\alpha(a-r)^2+\beta(a-r)^4+O\left((a-r)^6\right).
\end{equation*}
%
\bp\label{prop:phi-main}
The following asymptotic expansion holds 
\begin{equation}\label{eq:phi-expr-nu}
\varphi(\nu) = \frac{\pi}{\sqrt{2\alpha b}}+
\frac{\left(6 b\beta-\alpha^2\right)\pi}{8 \sqrt{2} b^{3/2}\alpha^{5/2}}(b-\nu)+O\left((b-\nu)^2\right),\qquad \text{for } \nu\to b^{-}.
\end{equation}
\ep
As an immediate corollary we recover the formula for the cut point along the equator.
\bc The  coordinates of the cut point along the equator satisfy $\thcut = \frac{\pi}{\sqrt{2 \alpha b}}$.
\ec
\begin{remark} The Gaussian curvature $K$ of a surface of revolution is constant on parallels, and is expressed as $K(r) = -\frac{m''(r)}{m(r)}$ (see \cite[p. 162]{do-carmo-curves}).
In particular $K(a) = -\frac{m''(a)}{m(a)} = \frac{2\alpha}{b}$ on the equator (from which one gets also $\alpha > 0$) and one can recover the first conjugate point along the equator via the Jacobi equation.
Indeed it is a general fact
that the first conjugate point along $\gamma$ is $\gamma(t_\mathrm{conj})$, where $t_\mathrm{conj}$
 is the first positive zero of the solution of the differential equation
\begin{equation*}
\ddot u(t) + K(\gamma(t))u(t) = 0,\qquad u(0) = 0, \dot u(0) = 1.
\end{equation*}
In our case $K$ is constant and the solution of this equation is given by
\begin{equation*}
u(t) = \sqrt{\frac{b}{2\alpha}}\sin\left(\sqrt{\frac{2\alpha}{b}}t\right),
\end{equation*}
whose first positive zero is $\tcut = t_\mathrm{conj} = \frac{\pi b}{\sqrt{2 b\alpha}}$.
On the equator $\vd t = b\ud\theta$, so we get indeed $\thcut = \frac{\tcut}{b} = \frac{\pi}{\sqrt{2 b\alpha}}$.\\
\end{remark}

The proof of  Proposition \ref{prop:phi-main} requires two elementary lemmas, whose proofs are postponed to Appendix \ref{a:lemmas}.
\bl\label{lem:power}
Let $U\subset \R^2$ be a neighborhood of $(0, 0)$ and let $f: U\to\R$ be a smooth function.
Let $f(x, y) = \sum_{i+j\leq n-1} a_{ij}x^i y^j + O((x^2+y^2)^{n/2})$ be its Taylor expansion.
Then for $y\geq 0$ the function defined by
\begin{equation*}
F(0)=\frac{\pi }{2}f(0,0),\qquad 
F(y) = \int_0^y\frac{f\left(x, y\right)\ud x}{\sqrt{y^2-x^2}},\quad\text{ for } y>0, \end{equation*}
is smooth and satisfies
\begin{equation}\label{eq:bigf}
F(y) = b_0+b_1 y+\cdots+b_{n-1}y^{n-1}+O(y^n),\qquad\text{where}\ b_k = \sum_{j=0}^k a_{j, k-j}\int_0^1 \frac{u^j\ud u}{\sqrt{1-u^2}}.
\end{equation}
\el

\bl\label{lem:f-smooth}
Let $g: (-\eps, \eps)\to\R$ be a smooth function. Let $f:(-\eps, \eps)^2\to\R$ be defined as
\begin{equation*}
f(x, y):=\left\{\begin{aligned}
&\frac{g(x)-g(y)}{x-y}\qquad&\text{if }x\neq y,\\
&g'(x)&\text{if }x=y.
\end{aligned}\right.
\end{equation*}
Then $f$ is a smooth function. If $g(x) = \sum_{k=0}^{n}a_k x^k + O(|x|^{n+1})$. Then
\begin{equation*}
f(x, y) = \sum_{k=1}^{n}a_k\sum_{i+j=k-1}x^i y^j + O((x^2+y^2)^{n/2}).
\end{equation*}
\el

\begin{proof}[Proof of Proposition \ref{prop:phi-main}]
Recall  that $\nu = m(R)$ and $R \to a^-$ as $\nu\to b^-$. Define $\wh{R}:= a-R$. We have
\begin{equation}\label{eq:phi-x-y-expression}
\begin{aligned}
\varphi(\nu) &= \int_R^a \frac{2m(R) \ud r}{m(r)\sqrt{m(r)^2-m(R)^2}} \\
&= \int_0^{\wh{R}} \frac{2m(a-\wh{R})\ud\wh{r}}{m(a-\wh{r})\sqrt{m(a-\wh{r})+m(a-\wh{R})}\sqrt{m(a-\wh{r})-m(a-\wh{R})}}.
\end{aligned}
\end{equation}
Being $b > 0$, it is clear that the map
\begin{equation*}
(\wh{r}, \wh{R})\mapsto \frac{2m(a-\wh{R})}{m(a-\wh{r})\sqrt{m(a-\wh{r})+m(a-\wh{R})}},
\end{equation*}
is a smooth function of two variables in a neighborhood of $(\wh{r}, \wh{R}) = 0$.
Its Taylor expansion can be easily computed from the Taylor expansion of $m$.
As for $(m(a-\wh{r})-m(a-\wh{R}))^{-1/2}$, it can be written as
\begin{equation*}
\frac{1}{\sqrt{\wh{R}^2-\wh{r}^2}}\sqrt{\frac{\wh{R}^2-\wh{r}^2}{\psi(\wh{r}^2)-\psi(\wh{R}^2)}},
\end{equation*}
(recall that $m(r) = \psi((a-r)^2)$).
By Lemma \ref{lem:f-smooth} we know that $\frac{\psi(\wh{R}^2)-\psi(\wh{r}^2)}{\wh{R}^2-\wh{r}^2}$ is a smooth function of $(\wh{r}^2, \wh{R}^2)$ and
\begin{equation*}
\frac{\psi(\wh{R}^2)-\psi(\wh{r}^2)}{\wh{R}^2-\wh{r}^2} = -\alpha+\beta(\wh{R}^2+\wh{r}^2)+O(\wh{R}^4+\wh{r}^4).
\end{equation*}
We have $\alpha > 0$, so $\sqrt{\frac{\wh{R}^2-\wh{r}^2}{\psi(\wh{r}^2)-\psi(\wh{R}^2)}}$ is also a smooth function
in a neighborhood of $(\wh{r}, \wh{R}) = 0$.

Summing up, the right hand side of \eqref{eq:phi-x-y-expression} has the form required in Lemma \ref{lem:power}.
Performing explicitly the long but straightforward computation described above we get the expression
\begin{equation*}
\phi(\nu) = \int_0^a\frac{f(\wh{r}, \wh{R})\ud \wh{r}}{\sqrt{\wh{R}^2-\wh{r}^2}},
\end{equation*}
with
\begin{equation*}
f(x, y) = \frac{\sqrt 2}{\sqrt{b\alpha}} + \frac{5 \alpha^2+2 b\beta}{2 \sqrt{2} (b\alpha)^{3/2}} x^2+
\frac{-3 \alpha^2+2 b\beta}{2 \sqrt{2} (b\alpha)^{3/2}} y^2 +O(x^4+y^4).
\end{equation*}

Using the fact that $\int_0^1\frac{\ud u}{\sqrt{1-u^2}} =
\frac{\pi}{2}$, and $\int_0^1\frac{u^2\ud u}{\sqrt{1-u^2}} = \frac{\pi}{4}$, we obtain the coefficients
\begin{equation*}
b_0 = \frac{\pi}{2}a_{00} = \frac{\pi}{\sqrt{2b\alpha}},\qquad \quad
b_{2} = \frac{\pi}{2}a_{02} + \frac{\pi}{4}a_{20} = \frac{\left(6 b\beta-\alpha^2\right) \pi }{8 \sqrt{2} (b\alpha)^{3/2}}.
\end{equation*}
In other words we have the expansion
\begin{equation}\label{eq:phi-expr-R}
\varphi(\nu) = \frac{\pi}{\sqrt{2 b\alpha}}+
\frac{\left(6 b\beta-\alpha^2\right) \pi }{8 \sqrt{2} (b\alpha)^{3/2}}(a-R)^2 + O\left((a-R)^4\right).
\end{equation}

It follows from the fact that $\alpha > 0$ 
that $(a-R)^2$ is a smooth function of $\nu = m(R)$ in a neighborhood of $\nu = b$. An explicit computation gives
\begin{equation*}
(a-R)^2=\frac{b-\nu}{\alpha }+\frac{\beta  (b-\nu)^2}{\alpha ^3}+O\left((b-\nu)^3\right),
\end{equation*}
which together with \eqref{eq:phi-expr-R} leads to \eqref{eq:phi-expr-nu}.
\end{proof}
\begin{remark}\label{r:31}
From this formula it is clear that the case $6b\beta = \alpha^2$ is going to be singular.
It is easy to compute that this is equivalent to $K''(a) = 0$, where $K(r_0)$ is the Gaussian curvature
on the parallel $r = r_0$.
\end{remark}

\subsection{Variations of the optimal geodesic}\label{s:main2}
Recall that the geodesic $c_b(t)$ follows the equator and reaches the cut-conjugate point
$(r, \theta) = (a, \thcut)$ for $t = \tcut$. We are now interested in variations of this optimal geodesic.
Such a variation is given by a smooth curve in the space of parameters,
that is a smooth curve $s \mapsto(\eta(s), t(s))\in\R^2$ such that $\eta(0) = 0$ and $t(0) = \tcut$.
Recall that $\eta(s) = \arccos\frac{\nu(s)}{b}$
is the angle between $\dot c_{\nu(s)}(0)$ and $\left(\frac{\partial}{\partial\theta}\right)_{(a, 0)}$.

\bp\label{prop:equator-variation}
There exist variations of the optimal geodesic $c_b$ such that
\bqn
\dist\left(c_{\nu(s)}(t(s)), (a, \thcut)\right) = O(s^3).
\eqnn
If the two-sphere is nonsingular (that is $K''(a)\neq 0$), then there exist no variation of the optimal geodesic $c_b$ such that
\bqn
\dist\left(c_{\nu(s)}(t(s)), (a, \thcut)\right) = o(s^3).
\eqnn
\ep
\bl\label{lem:degeneracy-main}
For $0 < \nu < b$ let $(r_\nu, \thcut)$ be the first point of intersection of $c_\nu$
with the meridian $P := \{\theta = \thcut\}$.
Then $r_\nu = a- \frac{\left(6 b\beta-\alpha^2\right)\sqrt b\pi}{16 \sqrt{2} \alpha^{5/2}}\eta^3+O(\eta^5)$.
\el
\begin{proof}
We assume that $c_\nu(t)$ reaches $r=R$ before it reaches $\theta=\thcut$ (this will be true for $\eta$ small enough). This means that $\dot r(t) \geq 0$ for $\theta \in [\thcut, \varphi(\nu)]$.
Thus from \eqref{eq:deriv-r} and the Clairaut relation we get
\begin{equation*}
\dot r(t) = \frac{\sqrt{m(r(t))^2-\nu^2}}{m(r(t))}, \qquad \quad 
\dot\theta(t) = \frac{\nu}{m(r(t))^2},
\end{equation*}
this permits to compute
\begin{align*}
\frac{\ud r}{\ud\theta} &= \frac{\dot r(t)}{\dot\theta(t)} = \frac{m(r(t))\sqrt{m(r(t))^2-\nu^2}}{\nu},\\
 \ddot r(t) &=\frac{\nu ^2 m'(r(t))}{m(r(t))^3},\\
 \ddot\theta(t) &=-\frac{2 \nu  \sqrt{m(r(t))^2 - \nu^2} m'(r(t))}{m(r(t))^4},  \\
\frac{\ud^2 r}{\ud\theta^2} &= \frac{\ddot r(t)\dot \theta(t)-\dot r(t)\ddot\theta(t)}{(\dot\theta(t))^3} = \frac{\left(2 m(r(t))^2-\nu^2\right) m(r(t)) m'(r(t))}{\nu ^2}.
\end{align*}
In particular
\begin{equation*}
\left(\frac{\ud r}{\ud\theta}\right)_{\theta=\varphi(\nu)} = \frac{b\sqrt{b^2-\nu^2}}{\nu} = b \tan(\eta),
\end{equation*}
and
\begin{equation*}
\frac{\ud^2 r}{\ud\theta^2} = O(R) = O(\eta).
\end{equation*}
Thus
\begin{align*}
a-r_\nu &= \int_{\thcut}^{\varphi(\nu)}\left(\frac{\ud r}{\ud\theta}\right)_{\theta_1}\vd\theta_1 =
\int_{\thcut}^{\varphi(\nu)}\left(\left(\frac{\ud r}{\ud\theta}\right)_{\varphi(\nu)}
-\int_{\theta_1}^{\varphi(\nu)}\left(\frac{\ud^2 r}{\ud\theta^2}\right)\vd\theta\right)\vd\theta_1\\
&= b\tan(\eta)(\varphi(\nu)-\thcut)+\int_{\thcut}^{\varphi(\nu)}\int_{\theta_1}^{\varphi(\nu)}O(\eta)\ud\theta\ud\theta_1.
\end{align*}
Recall that $b-\nu = b(1-\cos\eta) =\frac{b\eta^2}{2}+O(\eta^4)$.
The second term is $O(\eta^5)$ because we integrate twice on intervals of length $O(b-\nu) = O(\eta^2)$.
The conclusion follows by substituting $b-\nu = \frac{b\eta^2}{2}+O(\eta^4)$ in Proposition \ref{prop:phi-main}.
\end{proof}
%
\bl \label{lem:distance}
The distance between $(a, \thcut)$ and the geodesic segment
$c_\nu([0, \phi(\nu)])$ equals $\frac{\left(6 \beta b-\alpha^2\right)\sqrt b\pi}{16 \sqrt{2} \alpha^{5/2}}$\ $\eta^3+O(\eta^4)$.
\el
\begin{proof}
Let $\widetilde{q}$ be a point on the geodesic
segment under consideration.

Suppose that $\dist(\widetilde{q}, (a, \thcut)) = O(\eta^3)$.
Meridians are minimal curves, so we obtain $a - r_{\widetilde{q}} = O(\eta^3)$.
Hence, by the triangle inequality, $|\theta_{\widetilde{q}}-\thcut| = O(\eta^3)$. Repeating the computation from Lemma
\ref{lem:degeneracy-main} with $\theta_{\widetilde{q}}$ instead of $\thcut$ gives
\begin{equation*}
a-r_{\widetilde{q}} = \frac{\left(6 \beta b-\alpha^2\right)\sqrt b\pi}{16 \sqrt{2} \alpha^{5/2}}\eta^3+O(\eta^4).
\end{equation*}
Thus
\begin{equation*}
\dist(\widetilde{q}, (a, \thcut)) \geq \frac{\left(6 \beta b-\alpha^2\right)\sqrt b\pi}{16 \sqrt{2} \alpha^{5/2}}\eta^3+O(\eta^4).
\end{equation*}
\end{proof}
%
\begin{proof}[Proof of Proposition \ref{prop:equator-variation}]
The first statement is a direct consequence of Lemma \ref{lem:degeneracy-main}.

Now let $(\eta(s), t(s))$ be a variation of the equator. By the Gauss Lemma we get $t'(0) = 0$. Hence $\eta(s) \sim s$.
The conclusion follows from Lemma \ref{lem:distance}.
%
%
%
\end{proof}

Theorem \ref{thm:main} follows now from Proposition \ref{prop:equator-variation} and Corollary \ref{cor:the-exponent}.

\subsection{Oblate ellipsoid as a two-sphere of revolution}\label{sec:ellipsoid}
Let $M$ be an ellipsoid with semi-axes $b, b, c$, where $b \geq c$.
We denote $p, q$ its northern and southern pole respectively.
Let $a$ be the distance from the equator to a pole, so that $a$ and $b$ have the same meaning as before.
We will express the expansion of the function $r\mapsto m(r)$ and  $\nu\mapsto \phi(\nu)$
in terms of the semi-axes $b$, $c$. We will see that an oblate ellipsoid which is not a sphere is nonsingular.

\bp
Let M be an oblate ellipsoid with axes $b, b, c$, where $b > c$, is a two-sphere of revolution which satisfies assumptions (A1), (A2) and (A3).
The function $m$ has an expansion 
\bqn\label{eq:yyy}
m(r) = b - \alpha (a-r)^2+\beta (a-r)^4+O\left((a-r)^6\right)
\eqn
where $\alpha = \frac{b}{2c^2}$ and $\beta = \frac{b(4b^2-3c^2)}{24c^6}$.
\ep
\begin{proof} Clearly $M$ is a two-sphere of revolution satisfying assumptions (A1) and (A2).
It remains to prove \eqref{eq:yyy} and the fact $M$ is nonsingular.

Let $x$ be a point in the northern half of $M$ and let $\nu$ be its distance from the rotational axis.
Let $R$ be the geodesic distance from $x$ to $p$. Then $m(R) = \nu$.
On the other hand, $a-R$ is the length of an arc of an ellipse, which can be expressed by means of an elliptic integral.
In this way we obtain
\begin{equation*}
a-R = \int_\nu^b\sqrt{\frac{b^2-c^2}{b^2}+\frac{c^2}{b^2-t^2}}\ud t.
\end{equation*}
Let $Z = \sqrt{b-\nu}$. After substitution $t = b-z^2$ and some operations on power series the integral above transforms into
\begin{equation*}
\int_0^Z \frac{c\sqrt{2}}{\sqrt{b}}+\frac{\left(4 b^2-3 c^2\right) z^2}{2 \sqrt{2} b^{3/2} c}+O(z^4)\ud z,
\end{equation*}
which results in the expansion
\begin{equation*}
a-R = \frac{\sqrt{2} c}{\sqrt{b}}Z+\frac{\left(4 b^2-3 c^2\right)}{6 \sqrt{2} b^{3/2} c}Z^3 + O(Z^5).
\end{equation*}
By the Implicit Function Theorem $Z$ is a smooth function of $R$ in a neighborhood of $R = a$.
It follows from the expression above that
\begin{equation*}
\sqrt{b-\nu} = Z = \frac{\sqrt{b}}{\sqrt{2} c}(a-R)-\frac{\sqrt{b} \left(4 b^2-3 c^2\right)}{24 \sqrt{2} c^5}(a-R)^3 + O\left((a-R)^5\right),
\end{equation*}
so finally, squaring both sides,
\begin{equation*}
m(R) = \nu = b-\frac{b (a-R)^2}{2 c^2}+\frac{b \left(4 b^2-3 c^2\right)}{24 c^6}(a-R)^4 + O\left((a-R)^6\right).
\end{equation*}
This proves \eqref{eq:yyy}. 
Finally, notice that the singularity condition reads (see also Remark \ref{r:31})
\begin{equation*}
6b\beta-\alpha^2 = \frac{b^2(b^2-c^2)}{c^6}=0\qquad\text{if and only if}\qquad b = c.
\end{equation*}
Hence spheres are the only singular oblate ellipsoids of revolution. 
\end{proof}
\appendix
\section{Heat kernel asymptotics on $S^{2}$}\label{app:s2}
Let $x\in S^2$ be the north pole
and $\wh{x} \in S^2$ be the south pole of a standard sphere of radius 1.
After suitable renormalizations we obtain from formula (1) in \cite{fischer2s}
\bqn
p_t(x, \wh{x}) = \frac{1}{8\pi}\sum_{n=-\infty}^{\infty}(-1)^n(2n+1)e^{-(n^2+n)t}.
\eqn
Define $f(y)=\frac{e^{t/4}}{8t^{3/2}\pi^{1/2}}e^{-y^2/4t}iye^{-iy/2}$. Using standard rules of computing
Fourier transforms it is easy to see that
\bqn
\phi(n) := \frac{1}{2\pi}\int_{\R}f(y)e^{-iny}\ud y = \frac{1}{8\pi}(2n+1)e^{-(n^2+n)t}.
\eqnn
By Poisson Summation Formula we know that $\sum_{n=-\infty}^{\infty}\phi(n)e^{iny} = \sum_{k=-\infty}^{\infty}f(y+2k\pi)$. Using this identity for $y = \pi$ we obtain
\bqn
p_t(x, \wh{x}) = \sum_{n=-\infty}^{\infty}(-1)^n\phi(n) = \sum_{k=-\infty}^{\infty}f((2k+1)\pi).
\eqnn
Notice that for an asymptotic expansion only two terms in the sum on the right are significant. Hence we get
\bqn
p_t(x, \wh{x}) \simeq f(-\pi) + f(\pi) = \frac{e^{t/4}\sqrt{\pi}}{4t^{3/2}}e^{-\pi^2/4t}\sim \frac{1}{t^{3/2}} e^{-\dist^{2}(x,\wh{x})/4t}.
\eqnn
\section{Proof of Lemmas}\label{a:lemmas}
In this section we give a proof of Lemma \ref{lem:power} and Lemma \ref{lem:f-smooth}.

\begin{proof}[Proof of Lemma \ref{lem:power}]
The function $G(y)=\int_{0}^{1}\frac{f(ty,y)}{\sqrt{1-t^{2}}}dt$ is smooth on $(-\eps, \eps)$.
For $y > 0$ the change of variables $x = ty$ gives $F(y)=G(y)$. Moreover, by definition $F(0)=G(0)$,
hence $F$ is smooth on $[0,\eps)$.

We have $f(ty, y) = c_0(t) + c_1(t) y + \ldots + c_{n-1}(t)y^{n-1} + O(y^n)$, where $c_k(t) = \sum_{j=0}^k a_{j, k-j}t^j$.
Integrating with respect to $t$ we obtain $F(y) = G(y) = \wh{b}_0 + \wh{b}_1 y + \ldots + \wh{b}_{n-1}y^{n-1} + O(y^n)$,
where $\wh{b}_k = \int_0^1 \frac{c_k(t)\ud t}{\sqrt{1-t^2}}$, which is exactly the coefficient $b_k$ defined in the
statement of the Lemma.
\end{proof}

\begin{proof}[Proof of Lemma \ref{lem:f-smooth}]
Define $h: [0, 1]\times(-\eps, \eps)^2\to\R$ by $h(t, x, y) = g'(tx + (1-t)y)$.
Then $h(t, \cdot)$ is a continuous family of smooth functions. We have
\begin{equation*}
f(x, y) = \int_0^1 h(t, x, y)\ud t,
\end{equation*}
which proves that $f$ is smooth.
The Taylor expansion of $g$ is found by comparing the coefficients on both sides of the equality $(x-y)g(x, y) = f(x) - f(y)$.
\end{proof}

{\bf Acknowledgements.} The authors would like to thank Ugo Boscain for suggesting the problem and useful discussions.  This research has been supported by the European Research Council, ERC StG 2009 ``GeCoMethods'', contract number 239748, by the ANR Project GCM, program ``Blanche'', project number NT09-504490.

{\small
\bibliographystyle{siam}
\bibliography{bj}
}
\end{document}